\documentstyle[euscript,amsfonts,amssymb]{article}

\title{Meager-nowhere dense games (III): Remainder strategies.}

\author{Marion Scheepers\thanks{Supported in part by Idaho State Board of
   Education grant 91-093.}\\Department 
   of Mathematics,\\Boise State University,\\Boise, Idaho 83725}

\date{}

\newcommand{\QED}{\vrule width 6pt height 6pt depth 0pt \vspace{0.1in}}

\newcommand{\reals}{\Bbb R}

\newcommand{\naturals}{\Bbb N}

\newcommand{\wmg}{WMG(J)}
\newcommand{\wmeg}{WMEG(J)}

\newcommand{\vsgkl}{VSG([\kappa]^{<\lambda})}

\newcommand{\play}{(M_1,N_1,\dots,M_k,N_k,\dots)}

\newtheorem{theorem}{Theorem}
\newtheorem{problem}{Problem}
\newtheorem{prop}[theorem]{Proposition}
\newtheorem{lemma}[theorem]{Lemma}
\newtheorem{corollary}[theorem]{Corollary}

\begin{document}
\maketitle

\abstract{Player ONE chooses a meager set and player TWO, a nowhere
dense set per inning. They play
$\omega$ many innings. ONE's consecutive choices must form a
(weakly) increasing sequence. TWO wins if the union of the
chosen nowhere dense sets covers the union of the chosen meager sets.
A strategy of TWO which depends on knowing only the uncovered part of
the most recently chosen meager set is said to be a remainder strategy. 
TWO has a winning remainder strategy for
this game played on the real line with its usual topology.}

\section{Introduction}

   A variety of topological games from the class of meager-nowhere
   dense games were introduced in the papers \cite{B-J-S},
   \cite{S1} and \cite{S2}.  The existence
   of winning strategies which 
   use only the most recent move of either player (so-called coding
   strategies) and the existence of winning strategies which use only
   a bounded number of moves of the opponent as information (so-called
   $k$-tactics) are studied there and in
   \cite{K} and \cite{S3}. These studies are continued here
   for yet another fairly natural type of strategy, the so-called
   {\em remainder strategy}. 

   The texts \cite{E-H-M-R} and \cite{W} would be sufficient
   references for the miscellaneous results from combinatorial set theory
   which we use. As for notation: The symbol $J_{\reals}$ denotes the
   ideal of nowhere dense subsets of the real line (with its usual
   topology), while the symbol ``$\subset$'' is used exclusively to
   denote `` is a proper subset of ''. The symbol ``$\subseteq$'' is used
   to denote `` is a subset of, possibly equal to''. Let $(S,\tau)$ be
   a $T_1$-space 
   without isolated points, and let $J$ be its ideal of nowhere-dense subsets.
   The symbol $\langle J\rangle$ denotes the collection of
   meager subsets of the space. For $Y$ a subset of $S$, the symbol
   $J\lceil_Y$ denotes the set $\{T\in J:T\subseteq Y\}$.

   The game $\wmeg$ (defined in \cite{S2}) proceeds as follows: In the
   first inning player ONE chooses a meager set $M_1$, and player TWO responds
   with a nowhere dense set $N_1$. In the second inning player ONE chooses a
   meager set $M_2$, subject to the rule that $M_1\subseteq M_2$; TWO responds
   with a nowhere dense set $N_2$, and so on. The players play an inning for
   each positive integer, thus constructing a sequence $$\play$$
   which has the properties that $M_k\subseteq M_{k+1} \in\langle
   J\rangle$, and $N_k\in J$
   for each $k$. Such a sequence is said to be a {\em play of $\wmeg$}. Player
   TWO wins such a play if
\[\bigcup_{k=1}^{\infty}M_k=\bigcup_{k=1}^{\infty}N_k.\]
   A strategy of player TWO of the form 
\begin{enumerate}
\item{$N_1=F(M_1)$ and}
\item{$N_{k+1}=F(M_{k+1}\backslash(\bigcup_{j=1}^{k}N_j))$ for all $k$}
\end{enumerate}
   is said to be a {\em remainder strategy}. When does TWO have a winning
   remainder strategy in the game $\wmeg$?
  
   It is clear that TWO has a winning
   remainder strategy in $\wmeg$ if $J=\langle J\rangle$.
   The situation when $J\subset\langle
   J\rangle\subseteq{\EuScript P}(S)$ is not so easy. In Section 2 we
   investigate this question.
   We prove among other things Theorem
   \ref{wmeg-pos1}, which implies that TWO has a winning
   remainder strategy in the game $WMEG(J_{\reals})$.

   The game $\wmg$ proceeds just like $\wmeg$; only now the winning
   condition on TWO is relaxed so that TWO wins if
\[\bigcup_{n=1}^{\infty}M_n\subseteq\bigcup_{n=1}^{\infty}N_n.\]
   In Section 3 we study remainder strategies for this game; we
   briefly also discuss the game $SMG(J)$ here. In Section 4
   we attend to the version $VSG(J)$. The rules of this
   game turns out to be more advantageous to TWO from the point of view
   of existence of winning remainder strategies. 
   Some of the theorems in these two sections show that the hypotheses of
   Theorem \ref{wmeg-pos1} are to some extent necessary.

   Theorem \ref{countablefin-neg} is due to 
   Winfried Just, while Theorem \ref{countablefin2} is due to
   Fred Galvin. 
   I thank Professors Galvin and Just for kindly
   permitting me to present their result here and for fruitful
   conversations and correspondence concerning remainder strategies. 

\section{The weakly monotonic equal game, $WMEG(J)$.}

   When defining a remainder strategy $F$ for TWO in $WMEG(J)$, we shall
   take care that for each $A\in\langle J\rangle$:
\begin{enumerate}
\item{$F(A)\subseteq A$, and}
\item{$F(A)\neq\emptyset$ if (and only if) $A\neq\emptyset$.}
\end{enumerate}
   Otherwise, the strategy $F$ is sure not to be a winning remainder
   strategy for TWO in $WMEG(J)$,

\begin{theorem}\label{wmeg-pos1}If $(\forall X\in\langle
   J\rangle\backslash J)(cof(\langle J\rangle,\subset)\leq |J\lceil_X|)$,
   then TWO has a winning remainder strategy in $WMEG(J)$. 
\end{theorem}

   Theorem \ref{wmeg-pos1} follows from the next two lemmas.   In the proof of
   Lemma \ref{wmeg-pos1lemma} we use an
   auxilliary game, denoted 
   $REG(J)$. It is played as follows: A sequence $$\play$$ is a play of
   $REG(J)$ if $M_k\in\langle J\rangle$ and $N_k\in J$ for each $k$.
   Player TWO is declared the winner of a play of $REG(J)$ if
   $\cup_{k=1}^{\infty}M_k=\cup_{k=1}^{\infty}N_k$. TWO
   has a winning perfect information strategy in $REG(J)$. (We call
   $REG(J)$ the ``random equal game on $J$''.)

\begin{lemma}\label{wmeg-pos1lemma} If
\begin{enumerate}
\item{$cof(\langle J\rangle,\subset)$ is infinite and}
\item{$(\forall X\in\langle J\rangle\backslash J)(cof(\langle
   J\rangle,\subset)\leq |J\lceil_X|),$}
\end{enumerate}
   then TWO has a winning remainder strategy in $WMEG(J)$.
\end{lemma}

\begin{description}\item[Proof]{ Let ${\cal A}\subset\langle
   J\rangle\backslash J$ be a cofinal family of minimal cardinality.
   Observe that 
   $|{\cal A}|\leq|{\EuScript P}(X)|$ for each
   $X\in\langle J\rangle\backslash J$. Thus, if there is no $Y\in
   J\lceil_X$ such that $|{\cal A}|\leq|{\EuScript P}(Y)|$, then
   $|Y|<|X|$ for each $Y\in J\lceil_X$, and we fix a decomposition
\[X=\bigcup_{n=1}^{\infty}X_n\]
   where $\{X_n:n\in\naturals\}$ is a disjoint collection of sets from
   $\langle J\rangle\backslash J$. For each such $X_n$ we further fix a
   representation
\[X_n=\bigcup_{m=1}^{\infty}X_{n,m}\]
   where $X_{n,1}\subseteq X_{n,2}\subseteq\dots$ are from $J$, and a
   surjection 
\[\Theta^X_n:J\lceil_{X_n}\rightarrow\mbox{$^{<\omega}{\cal A}$}.\]

   For each $Y\in J$ such that $|{\cal A}|\leq|{\EuScript P}(Y)|$ the
   set $Y$ is infinite and we also write
\[Y=\bigcup_{n=1}^{\infty}Y_n\]
   where $\{Y_n:n\in\naturals\}$ is a pairwise disjoint collection
   such that $|Y_n|=|Y|$ for each $n$. Further, choose for each $n$ a
   surjection 
\[\Psi^Y_n:{\EuScript P}(Y_n)\backslash\{\emptyset,Y_n\}\rightarrow
   \mbox{$^{<\omega}{\cal A}$}.\]

   Let $U$ and $V$ be sets in $\langle J\rangle$ such that we have
   chosen a decomposition $U=\cup_{n=1}^{\infty}U_n$ as above. We'll use
   the notation
\[U\subseteq^*V\]
   to denote that there is an $m$ such that $U_n\subseteq V$ for each
   $n\geq m$; we say that $m$ {\em witnesses that} $U\subseteq^*V$.

   Fix a well-ordering $\prec$ of $\langle J\rangle$. For $X\in\langle
   J\rangle$ we define:
\begin{enumerate}
\item{$\Theta(X)$: the $\prec$-minimal element $A$ of ${\cal A}$ such that
   $X\subseteq A$,}
\item{$\Phi(X)$: the $\prec$-minimal element $Z$ of $\langle
   J\rangle\backslash J$ such that $Z\subseteq^*X$ whenever this is
   defined, and the empty set otherwise,}
\item{$k(X)$: the minimal natural number which witnesses that
   $\Phi(X)\subseteq^*X$ whenever $\Phi(X)\neq\emptyset$, and $0$ otherwise,}
\item{$\Gamma(X)$: the $\prec$-minimal $Y\in J$ such that
   $|J\lceil_X|\leq |{\EuScript P}(Y)|$ and $Y\subseteq^*X$ whenever this
   is defined, and the empty set otherwise, and}
\item{$m(X)$: the minimal natural number which witnesses that
   $\Gamma(X)\subseteq^*X$ whenever $\Gamma(X)\neq\emptyset$, and $0$
   otherwise.}
\end{enumerate}
   Let $G$ be a winning perfect information strategy for TWO in the
   game $REG(J)$. We are now ready to define TWO's remainder strategy
\[F:\langle J\rangle\rightarrow J.\]
   Let $B\in\langle J\rangle$ be given.\\
\underline{$B\in J$:} Then we define $F(B)=B$.\\
\underline{$B\not\in J$:} We distinguish between two cases:\\
{\bf Case 1:} $\Gamma(B)=\emptyset$.\\
   Write $X$ for $\Phi(B)$ and $n+1$ for $k(B)$. For $1\leq j\leq n$
   define $\sigma_j$ so that
\[\sigma_j=\left\{\begin{array}{ll}
   \Theta^X_j(X_j\backslash B) & \mbox{if $X_j\backslash B\in J$}\\
   \emptyset & \mbox{otherwise}
   \end{array}\right.
\]
   Let $\tau$ be
   $\sigma_1\frown\dots\frown\sigma_n\frown\langle\Theta(B)\rangle$,
   the concatenation of these finite sequences, and choose $V\in
   J\lceil_{X_{n+1}}$ such that $\Theta^X_{n+1}(V)=\tau$. Then define
\[F(B)=B\cap[X_{1,n+1}\cup\dots\cup X_{n,n+1}\cup V\cup
   ((\cup\{G(\sigma):\sigma\subseteq \tau\})\backslash X)].
\]
{\bf Case 2:} $\Gamma(B)\neq\emptyset$.\\
   Write $Y$ for $\Gamma(B)$ and $n+1$ for $m(B)$. For $1\leq j\leq n$
   define $\sigma_j$ so that
\[\sigma_j=\left\{\begin{array}{ll}
   \Psi^Y_j(Y_j\backslash B) & \mbox{if $Y_j\backslash
   B\not\in\{\emptyset,Y_j\}$}\\ 
   \emptyset & \mbox{otherwise}
   \end{array}\right.
\]
   Let $\tau$ be
   $\sigma_1\frown\dots\frown\sigma_n\frown\langle\Theta(B)\rangle$, 
   the concatenation of these finite sequences, and choose
   $V\in{\EuScript P}(Y_{n+1})\backslash\{\emptyset, Y_{n+1}\}$ so that
   $\Psi^Y_{n+1}(V)=\tau$. Then define
\[F(B)=B\cap[Y_1\cup\dots\cup Y_n\cup V\cup
   ((\cup\{G(\sigma):\sigma\subseteq \tau\})\backslash Y)].
\]

   From its definition it is clear that $F(B)\subseteq B$ for each
   $B\in\langle J\rangle$. To see that $F$ is a winning remainder
   strategy for TWO in $WMEG(J)$, consider a play 
$$\play$$ 
   during which TWO adhered to the strategy $F$. To facilitate the
   exposition we write:
\begin{enumerate}
\item{$B_1$ for $M_1$ and $B_{j+1}$ for
   $M_{j+1}\backslash\cup_{i=1}^jN_i$,}
\item{$Y^j$ for $\Gamma(B_j)$,}
\item{$X^j$ for $\Phi(B_j)$,}
\item{$A^j$ for $\Theta(B_j)$,}
\item{$k_j$ for $k(B_j)$ and }
\item{$m_j$ for $m(B_j)$.}
\end{enumerate}

   We must show that
   $\cup_{j=1}^{\infty}B_j\subseteq\cup_{j=1}^{\infty}N_j$. 
   
   We may assume that $B_j\not\in J$ for each $j$.

   Suppose that $Y^{j+1}\neq\emptyset$ for some $j$. Then $N_{j+1}$ is
   defined by Case 2, and as such is of the form
\[B_{j+1}\cap[Y^{j+1}_1\cup\dots\cup Y^{j+1}_{m_{j+1}-1}\cup
V_{j+1}\cup((\cup\{ G(\sigma):\sigma\subseteq \tau_{j+1}\})\backslash Y^{j+1})]
\]
   where $V_{j+1}$ and $\tau_{j+1}$ have the obvious meanings. Now
   $B_{j+1}\backslash B_{j+2}=N_{j+1}$, and thus:
\begin{itemize}
\item{$Y^{j+1}\subseteq^*B_{j+2}$ is a candidate for $Y^{j+2}$, and}
\item{$Y^{j+2}\neq\emptyset$, so that $N_{j+2}$ is also defined by Case 2.}
\end{itemize}
   We conclude that if $Y^j\neq\emptyset$ for some $j$, then
   $Y^i\neq\emptyset$ for all $i\geq j$. Moreover, it then also
   follows that $Y^{i+1}\preceq Y^i$ for each $i\geq j$. Since $\prec$
   is a well-order, there
   is a fixed $k$ such that $Y^i=Y^k$ for all $i\geq k$. Let $Y$ be
   this common value of $Y^i, i\geq k$. An inductive computation now
   shows that $(A^k,\dots,A^j)\subseteq\tau_j$ for each $j\geq k$. But
   then
\[B_j\cap[(G(A^k)\cup\dots\cup G(A^k,\dots,A^j))\backslash Y]\subseteq N_j\]
   for each $j\geq k$. It follows that $\cup_{j=k}^{\infty}B_j\backslash
   Y\subseteq \cup_{j=k}^{\infty}N_j$. But it is also clear that
   $Y\cap(\cup_{j=1}^{\infty}B_j)\subseteq\cup_{j=k}^{\infty}N_j$. It
   then follows from the monotonicity of the sequence of $M_j$-s that TWO
   has won this play.
   
   The other case to consider is that $Y^{j+1}=\emptyset$ for all $j$.
   In this case, $X^{j+1}\neq\emptyset$ for each $j$, and $N_{j+1}$ is
   defined by Case 1. In this case $N_{j+1}$ is of the form:
\[B_{j+1}\cap[X^{j+1}_{1,k_{j+1}}\cup\dots\cup
X^{j+1}_{k_{j+1}-1,k_{j+1}} \cup V_{j+1}
\cup((\cup\{G(\sigma):\sigma\subseteq\tau_{j+1}\})\backslash X^{j+1})],\]
   where $V_{j+1}$ and $\tau_{j+1}$ have the obvious meaning.
   Now $B_{j+1}\backslash B_{j+2}=N_{j+1}$, and so
   $X^{j+1}\subseteq^*B_{j+2}$, and $X^{j+1}$ is a candidate for $X^{j+2}$.
   It follows that $X^{j+2}\preceq X^{j+1}$ for each $j<\omega$.
   Since $\prec$ is a well-order we once again fix $k$ such that
   $X^j=X^k$ for all $j\geq k$. Let $X$ denote $X^k$. As before, $\langle
   A^k,\dots,A^j\rangle\subseteq\tau_j$ for each such $j$, and it follows
   that TWO also won these plays. $\QED$}
\end{description}

\begin{lemma}\label{wmeg-pos2lemma} If $\langle J\rangle={\EuScript
   P}(S)$, then TWO has a winning remainder strategy in $WMEG(J)$.
\end{lemma}

\begin{description}\item[Proof]{ Let $\prec$ be a well-order of
   ${\EuScript P}(S)$, and write $S=\cup_{n=1}^{\infty}S_n$ such that
   $S_n\in J$ for each $n$, and the $S_n$-s are pairwise disjoint. For
   each countably infinite $Y\in J$ write 
   $Y=\cup_{n=1}^{\infty}Y_n$ so that $|Y_n|=n$ for each $n$, and
   $\{Y_n:n\in\naturals\}$ is pairwise disjoint. For $X$ and $Y$ in
   $\langle J\rangle$ write $Y\subseteq^*X$ if $Y\backslash X$ is finite.

   For each $X\in\langle J\rangle\backslash J$,
\begin{itemize}
\item{either there is an infinite $Y\in J\lceil_X$,}
\item{or else $X$ is countably infinite.}
\end{itemize}

   In the first of these cases, let $\Phi(X)$ be the $\prec$-first
   countable element $Y$ of $J$ such that $Y\subseteq^*X$, and let $m(X)$
   be the smallest $n$ such that $Y_m\subseteq X$ for all $m\geq n$. 

   In the second of these cases, let $\Phi(X)$ be the $\prec$-least element
   $Y$ of $\langle J\rangle\backslash J$ such that $Y\subseteq^*X$, and
   let $m(X)$ be the minimal $n$ such that $\Phi(X)\cap S_m\subseteq X$
   for all $m\geq n$. Also write $\Phi(X)_j$ for $\Phi(X)\cap S_j$ for
   each $j$, in this case.

   Then define $F(X)$ so that
\begin{enumerate}
\item{$F(X)=X$ if $X\in J$, and}
\item{$F(X)=X\cap[(S_1\cup\dots\cup
   S_{m(X)})\backslash\Phi(X))\cup(\Phi(X)_1 \cup\dots\cup\Phi(X)_{m(X)}]$}
\end{enumerate}

   Then $F$ is a winning remainder strategy for TWO in $WMEG(J)$, for
   reasons analogous to those in the proof of Lemma \ref{wmeg-pos1lemma}.
   $\QED$}
\end{description}

\begin{corollary}\label{cor1} Player TWO has a winning remainder
   strategy in the game $WMEG(J_{\reals})$.
\end{corollary}

   Recall (from \cite{S2}) that $G$ is a coding strategy for
   TWO if:
\begin{enumerate}
\item{$N_1=G(\emptyset,M_1)$ and}
\item{$N_{k+1}=G(N_k,M_{k+1})$ for each $k$.}
\end{enumerate}

   If $F$ is a winning remainder strategy for TWO in $WMEG(J)$,
   then the function $G$ which is defined so that $G(W,B)=W\cup
   F(B\backslash W)$ is a winning coding strategy for TWO in $WMEG(J)$.
   Thus, Corollary \ref{cor1} solves Problem 2 of \cite{S2} positively.

   We shall later see that the sufficient condition for the
   existence of a winning coding strategy given in Theorem
   \ref{wmeg-pos1} is to some extent necessary (Theorems
   \ref{countablefin-neg} and \ref{countablefin2}). However, this condition
   is not absolutely necessary, as we shall now illustrate. First, note
   that for any decomposition $S=\cup_{j=1}^kS_k$, the following
   statements are equivalent:
\begin{enumerate}
\item{TWO has a winning remainder strategy in $WMEG(J)$,}
\item{For each $j$, TWO has a winning remainder strategy in
   $WMEG(J\lceil_{S_j})$.}
\end{enumerate}

   Now let $S$ be the disjoint union of the real line and a countable set
   $S^*$. Let $X\in J$ if $X\cap S^*$ is finite and $X\cap \reals\in
   J_{\reals}$. Then $S^*\in\langle J\rangle$, and $J\lceil_{S^*}$ is a
   countable set, while $cof(\langle J\rangle,\subset)$ is uncountable.
   According to Corollary \ref{cor1} and Lemma \ref{wmeg-pos2lemma},
   TWO has a winning remainder strategy in $WMEG(J)$.

   Let $\lambda$ be an infinite cardinal of countable cofinality. For
   $\kappa\geq \lambda$, we declare a subset of $\kappa$ to be open if it
   is either empty, or else has a complement of cardinality less than
   $\lambda$. With this topology, $J=[\kappa]^{<\lambda}$.

\begin{corollary}\label{cor2} Let $\lambda$ be a cardinal
   of countable cofinality, and let $\kappa>\lambda$ be a cardinal
   number. If $cof([\kappa]^{\lambda},\subset)\leq
   \lambda^{<\lambda}$, then TWO has a winning 
   remainder strategy in $WMEG([\kappa]^{<\lambda})$. 
\end{corollary}

   Let ${\cal A}$ be a subset of $\langle J\rangle$. The game
   $WMEG({\cal A},J)$ is played like $\wmeg$, except 
   that ONE is confined to choosing meager sets which are in ${\cal
   A}$ only. Thus, $WMEG(J)$ is the special case of $WMEG({\cal A},J)$ for
   which ${\cal A}=\langle J\rangle$. If there is a cofinal family
   ${\cal A}\subset \langle J\rangle$ such that TWO does not have a
   winning remainder strategy in $WMEG({\cal A},J)$, then TWO does not 
   have a winning remainder strategy in $WMEG(J)$.  

\begin{theorem}There is a cofinal family ${\cal
   A}\subset[\omega_1]^{\aleph_0}$ such that TWO does not have a winning
   remainder strategy in 
   $WMEG({\cal A},[\omega_1]^{<\aleph_0})$.
\end{theorem}

\begin{description}\item[Proof]{ Let $F$ be a remainder strategy for
   TWO such that $F(X)\subseteq X$ for every countable subset $X$ of
   $\omega_1$.
   Put ${\cal A}=\{\alpha<\omega_1: cof(\alpha)=\omega\}$.

   We show that there is a sequence 
$$\langle (S_1,T_1), (S_2,T_2),\dots\rangle$$
   such that:
\begin{enumerate}
\item{$S_n$ is a stationary subset of $\omega_1$,}
\item{$T_n$ is a finite subset of $\omega_1$,}
\item{$S_{n+1}\subseteq S_n$ for each $n$,}
\item{$F(\gamma)=T_1$ for each $\gamma\in S_1$ and}
\item{$F(\gamma\backslash(\cup_{j=1}^nT_j))=T_{n+1}$ for each
   $\gamma\in S_{n+1}$ and for each $n$.}
\end{enumerate}

   To establish the existence of $S_1$ and $T_1$ we argue as follows.

   For each $\gamma<\omega_1$ which is of countable cofinality we put
$$\phi_1(\gamma)=\max F(\gamma) ( <\gamma ).$$
   By Fodor's lemma there is a stationary set $S_0$ of countable limit
   ordinals, and an ordinal $\delta_0<\min(S_0)$ such that
   $\phi_1(\gamma)=\delta_0$ for each $\gamma\in S_0$. But then
   $F(\gamma)$ is a finite subset of $\delta_0+1$ for each such $\gamma$.
   Since every partition of a stationary set into countably many sets has
   at least one of these sets stationary, we find a stationary set
   $S_\subset S_0$ and a finite set $T_1\subset \delta_0+1$ such that
   $F(\gamma)=T_1$ for each $\gamma\in S_1$.

   This specifies $(S_1,T_1)$. Now let $1\leq n<\omega$ be given and
   suppose that $(S_1,T_1),\dots,(S_n,T_n)$ with properties 1 through 5
   are given.

   For $\gamma\in S_n$ we define:
$$\phi_{n+1}(\gamma)=\max F(\gamma\backslash(T_1\cup\dots\cup T_n))
(<\gamma).$$ 
   Once again there is, by Fodor's Lemma, a stationary set $S'\subset
   S_n$ and an ordinal $\delta'<\omega-1$ such that
   $\phi_{n+1}(\gamma)=\delta'$ for each $\gamma\in S'$. As before we
   then find a stationary set $S_{n+1}\subseteq S'$ and a finite set
   $T_{n+1}\subseteq \delta'+1$ such that
   $F(\gamma\backslash(T_1\cup\dots\cup T_n))=T_{n+1}$ for each
   $\gamma\in S_{n+1}$. 

   Then $(S_1,T_1),\dots, (S_{n+1},T_{n+1})$ have properties 1 through 5.
   It follows that there is an infinite sequence of the required sort. 

   Put $\delta=\sup(\cup_{n=1}^\infty T_n)$ and choose $\gamma_n\in S_n$
   such that 
$$\delta<\gamma_1<\gamma_2<\gamma_3<\dots$$
   Then $(\gamma_1, T_1, \gamma_2, T_2, \dots)$ is an $F$-play of
   $WMEG({\cal A},[\omega_1]^{<\aleph_0})$, and is lost by TWO. $\QED$}
\end{description}

   Though there may be cofinal families ${\cal A}$ such that TWO does
   not have a winning remainder strategy in $WMEG({\cal A},J)$, there may
   for this very same $J$ also be cofinal families ${\cal
   B}\subset\langle J\rangle$ such that TWO does have a winning remainder
   strategy in $WMEG({\cal B},J)$.

\begin{theorem}\label{wmeg-pos2} Let $\lambda$ be an infinite cardinal
   number of countable cofinality. If $\kappa>\lambda$ is a cardinal
   for which $cof([\kappa]^{\lambda},\subset)=\kappa$, then there is a
   cofinal family ${\cal A}\subset[\kappa]^{\lambda}$ such that TWO has a
   winning remainder strategy in $WMEG({\cal A},[\kappa]^{<\lambda})$.
\end{theorem}

\begin{description}\item[Proof]{ Let $(B_{\alpha}:\alpha<\kappa)$
   bijectively enumerate a cofinal subfamily of $[\kappa]^{\lambda}$.
   Write 
$$\kappa=\cup_{\alpha<\kappa}S_{\alpha}$$
   where $\{S_{\alpha}:\alpha<\kappa\}\subset[\kappa]^{\lambda}$ is a
   pairwise disjoint family.
 
   Define: $A_{\alpha}=\{\alpha\}\cup(\cup_{x\in B_{\alpha}}S_x)$ for each
   $\alpha<\kappa$, and put ${\cal A}=\{A_{\alpha}:\alpha<\kappa\}$. Then
   ${\cal A}$ is a cofinal subset of $[\kappa]^{\lambda}$. Also let
   $\Psi:{\cal A}\rightarrow\kappa$ be such that
   $\Psi(A_{\alpha})=\alpha$ for each $\alpha\in\kappa$.

   Choose a sequence
   $\lambda_1<\lambda_2<\dots<\lambda_n<\dots$ of cardinal numbers
   converging to $\lambda$. For each $A\in{\cal A}$ we write
   $A=\cup_{n=1}^{\infty}A^n$ where $A^1\subset A^2\subset\dots$ are such
   that $|A^n|=\lambda_n$ for each $n$.

   For convenience we write, for $C$ and
   $D$ elements of $[\kappa]^{\leq\lambda}$, that $C=^*D$ if $|C\Delta
   D|<\lambda$. Observe that for $A$ and $B$ elements of ${\cal A}$,
   $A\neq B$ if, and only if, $|A\Delta B|=\lambda$. 

   Now define TWO's remainder strategy $F$ as follows:

\begin{enumerate}
\item{$F(A)=\{\Psi(A)\}\cup A^1$ for $A\in{\cal A}$,}
\item{$F(A)=\{\Psi(B)\}\cup(\cup(\{C^{m+1}:\Psi(C)\in\Gamma(A)\})\cap
   B)\cup B^{m+1}$ if $A\not\in{\cal A}$ but $A\subset B$ and $A=^*B$ for
   some $B\in{\cal A}$. Observe that this $B$ is unique. In this
   definition, $\Gamma(A)=B\backslash A $, and $m$ is minimal such that
   $|\Gamma(A)|\leq \lambda_m$.}
\item{$F(A)=\emptyset$ in all other cases.}
\end{enumerate}

   Observe that $|F(A)|<\lambda$ for each $A$, so that $F$ is a
   legitimate strategy for TWO. To see that $F$ is indeed a winning
   remainder strategy for TWO, consider a play
$$(M_1,N_1,\dots,M_k,N_k,\dots)$$
   of $WMEG({\cal A},[\kappa]^{<\lambda})$ during which TWO used $F$
   as a remainder strategy.

   Write $M_i=A_{\alpha_i}$ for each $i$. By the rules of the game we
   have:
$$A_{\alpha_1}\subseteq A_{\alpha_2}\subseteq\dots.$$   
   Also, $N_1=\{\alpha_1\}\cup A_{\alpha_1}^1$ and $n_1$ minimal is such that
   $|N_1|\leq\lambda_{n_1}$. An inductive computation shows that 
   $N_{k+1}=F(M_{k+1}\backslash(\cup_{j=1}^k N_j)$ is the set
$$([\{\alpha_{k+1}\}\cup(\cup\{A_{\gamma}^{n_{k}+1}:\gamma\in
   N_k\})\cup A_{\alpha_{k+1}}^{n_k+1})\cap A_{\alpha+{k+1}}$$ 
   from which it follows that:
\begin{enumerate}
\item{$N_1\subseteq N_2\subseteq\dots\subseteq N_k\subseteq\dots$,}
\item{$n_1\leq n_2\leq\dots\leq n_k\leq\dots$ goes to infinity,}
\item{$\alpha_j\in N_k$ whenever $j\leq k$, and thus}
\item{$A_{\alpha_j}^p\subseteq N_k$ for $j\leq k$ and $p\leq
   n_{k-1}$.}
\end{enumerate}
   The result follows from these remarks. $\QED$}
\end{description}

   Theorem \ref{wmeg-pos2} also covers the case when
   $\lambda=\aleph_0$. For cofinal families ${\cal A}\subset\langle
   J\rangle$ which have the special property that 
$$A\neq B\Leftrightarrow A\Delta B\not\in J$$ 
   (like the one exhibited in the above proof), there is indeed an
   equivalence between the existence of winning coding strategies and
   winning remainder strategies in the game $WMEG({\cal A},J)$.
   Particularly:

\begin{prop}\label{wmeg-equiv} Let ${\cal A}\subset \langle
   J\rangle$ be a cofinal family such that for $A$ and $B$ elements of
   ${\cal A}$, $A\neq B\Leftrightarrow A\Delta B\not\in J$. Then the
   following statements are equivalent:
\begin{enumerate}
\item{TWO has a winning coding strategy in $WMEG({\cal A},J)$.}
\item{TWO has a winning remainder strategy in $WMEG({\cal A},J)$.}
\end{enumerate}
\end{prop}
\begin{description}\item[Proof]{We must verify that $1$ implies $2$.
   Thus, let $F$ be a winning coding strategy for TWO in the game
   $WMEG({\cal A},J)$. We define a remainder strategy $G$. Let $X$ be
   given. If $X\in{\cal A}$ we define $G(X)=F(\emptyset,X)$. If
   $X\not\in{\cal A}$ but there is an $A\in{\cal A}$ such that $X\subset
   A$ and $X=^*A$, then by the property of ${\cal A}$ there is a unique
   such $A$ and we set $T=A\backslash X(\in J)$. In this case define
   $G(X)=F(T,A)$. In all other cases we put $G(X)=\emptyset$. Then $G$ is
   a winning remainder strategy for TWO in $WMEG({\cal A},J)$. $\QED$}
\end{description}

   It is not always the case that there is a cofinal ${\cal
   A}\subset\langle J\rangle$ which satisfies the hypothesis of
   Proposition \ref{wmeg-equiv}. For example, let $J\subset{\EuScript
   P}(\omega_2)$ be defined so that $X\in J$ if, and only if,
   $X\cap\omega$ is finite and $X\cap(\omega_2\backslash\omega)$ has
   cardinality at most $\aleph_1$. Let
   $\{S_{\alpha}:\alpha<\omega_2\}$ be a cofinal family. Choose
   $\alpha\neq\beta\in\omega_2$ such that:
\begin{enumerate}
\item{$\omega\subset (S_{\alpha}\cap S_{\beta})$ and}
\item{$S_{\alpha}\neq S_{\beta}$.}
\end{enumerate}
   Then $S_{\alpha}\Delta S_{\beta}\in J$.

   Coupled with Theorem \ref{wmeg-pos2} and an assumption about
   cardinal arithmetic, the following Lemma (left to the reader)
   enables us to conclude much more.  

\begin{lemma}\label{inductlemma} Let $\lambda$ be a cardinal of
   countable cofinality, let $\mu\leq\lambda$ be a regular cardinal
   number, and let $\{S_{\alpha}:\alpha<\mu\}$ be a collection of
   pairwise disjoint sets such that for each $\alpha<\mu$ there is a
   cofinal family ${\cal A}_{\alpha}\subset[S_{\alpha}]^{\lambda}$ for
   which TWO has a winning remainder strategy in the game $WMEG({\cal
   A}_{\alpha}, [S_{\alpha}]^{<\lambda})$. Then there is a cofinal family
   ${\cal A}\subset[\cup_{\alpha<\mu}S_{\alpha}]^{\lambda}$ such that TWO
   has a winning remainder strategy in $WMEG({\cal
   A},[\cup_{\alpha<\mu} S_{\alpha}]^{<\lambda})$.
\end{lemma}

\begin{corollary}\label{wmegcor} Assume the Generalized Continuum
   Hypothesis. Let 
   $\lambda$ be a cardinal of countable cofinality. For every infinite
   set $S$ there is a cofinal family ${\cal A}\subset [S]^{\lambda}$ such
   that TWO has a winning remainder strategy in $WMEG({\cal
   A},[S]^{<\lambda})$.
\end{corollary}

    It is clear that if TWO has a winning remainder strategy in the
    game $WMEG(J)$, then TWO has a winning remainder strategy in $WMG(J)$.
    The converse of this assertion is not so clear.

\begin{problem} Is it true that if TWO has a winning remainder
   strategy in the game $WMG(J)$, then TWO has a winning remainder
   strategy in $WMEG(J)$?
\end{problem}

\section{The strongly monotonic game, $SMG(J)$.}

   A sequence $\play$ is a play of the strongly monotonic game if:
\begin{enumerate}
\item{$M_k\cup N_k\subseteq M_{k+1}\in\langle J\rangle$, and}
\item{$N_k\in J$ for each $k$.}
\end{enumerate}
   Player TWO wins such a play if
   $\cup_{j=1}^{\infty}M_j=\cup_{j=1}^{\infty}N_j$. This game was studied
   in \cite{B-J-S} and \cite{S1}; from the point of view of TWO this
   gives TWO a little more control over how ONE's meager sets increase as
   the game progresses. It is clear that if TWO has a winning remainder
   strategy in $WMG(J)$, then TWO has a winning remainder strategy in
   $SMG(J)$. The converse is also true, showing that in the context of
   remainder strategies, the more stringent requirements placed on ONE by
   the rules of the strongly monotonic game is not of any additional
   strategic value for TWO:

\begin{lemma}\label{wmgsmgequiv}If TWO has a winning remainder
   strategy in $SMG(J)$, then TWO has a winning remainder strategy in
   $WMG(J)$.
\end{lemma}
\begin{description}\item[Proof]{ Let $F$ be a winning remainder
   strategy for TWO in $SMG(J)$. We show that it is also a winning
   remainder strategy for TWO in $WMG(J)$.

   Let $\play$ be a play of $WMG(J)$ during which TWO used $F$ as a
   remainder strategy. Put $M^*_1=M_1$ and
   $M^*_{k+1}=M_{k+1}\cup(N_1\cup\dots\cup N_k)$ for each $k$. Then
   $(M^*_1,N_1,\dots,M^*_k,N_k,\dots)$ is a play of $SMG(J)$ during which
   TWO used the winning remainder strategy $F$. It follows that 
   $\cup_{k=1}^{\infty}M_k\subseteq\cup_{k=1}^{\infty}N_k$, so that TWO won
   the $F$-play of $WMG(J)$. $\QED$}
\end{description}

   We now restrict ourselves to the rules of $WMG(J)$.
   As with $WMEG(J)$, a winning remainder strategy for TWO in the game
   $WMG(J)$ gives rise to the existence of a winning coding strategy for
   TWO. In general, the statement that player TWO has a
   winning remainder strategy in the game $WMG(J)$ is stronger than the
   statement that TWO has a winning coding strategy.
   To see this, recall that TWO has a winning coding strategy in
   $WMG([\omega_1]^{<\aleph_0})$ (see Theorem 2 of \cite{S2}). But
   according to the next theorem, TWO does not have a winning
   remainder strategy in the game $WMG([\omega_1]^{<\aleph_0})$.

\begin{theorem}[Just]\label{countablefin-neg} If $\kappa\geq\aleph_1$,
   then TWO does not have a winning remainder strategy in the game
   $WMG([\kappa]^{<\aleph_0})$. 
\end{theorem}

\begin{description}\item[Proof]{ Let $F$ be a remainder strategy for
   TWO. For each $\alpha<\omega_1$ we put
$$\Phi(\alpha)=\sup(\cup\{F(\alpha\backslash T):
   T\in[\alpha]^{<\aleph_0}\}\cup\alpha).
$$ 
   Then $\Phi(\alpha)\geq\alpha$ for each such $\alpha$. Choose a closed,
   unbounded set $C\subset\omega_1$ such that:
\begin{enumerate}
\item{$\Phi(\gamma)<\alpha$ whenever $\gamma<\alpha$ are elements of
   $C$, and}
\item{each element of $C$ is a limit ordinal.}
\end{enumerate}

   Then, by repeated use of Fodor's pressing down lemma,  we inductively
   define a sequence $((\phi_1,S_1,T_1),\dots,(\phi_n,S_n,T_n),\dots)$
   such that: 
\begin{enumerate}
\item{$C\supset S_1\supset\dots\supset S_n\supset\dots$ are stationary
   subsets of $\omega_1$,}
\item{$F(\alpha)\cap\alpha=T_1$ for each $\alpha\in S_1$, and}
\item{$F(\alpha\backslash(T_1\cup\dots\cup T_n))=T_{n+1}$ for each $n$
   and each $\alpha\in S_n$.}
\end{enumerate}

   Put $\xi=\sup(\cup_{n=1}^{\infty}T_n)+\omega$. Choose $\alpha_n\in
   S_n$ so that $\xi\leq\alpha_1<\alpha_2<\dots<\alpha_n<\dots$. By the
   construction above we have:
\begin{enumerate}
\item{$F(\alpha_1)\cap\xi=T_1$ and}
\item{$F(\alpha_{n+1}\backslash(T_1\cup\dots\cup T_n))\cap\xi=T_{n+1}$
   for each $n$.}
\end{enumerate}

   But then
$(\cup_{n=1}^{\infty}T_n)\cap\xi\subset\xi=(\cup_{n=1}^{\infty}\alpha_n)\cap\xi$,
   so that TWO lost this play of $WMG([\omega_1]^{<\aleph_0})$. $\QED$}
\end{description}

   For a cofinal family ${\cal A}\subseteq\langle J\rangle$, the game
   $WMG({\cal A},J)$ proceeds just like $WMG(J)$, except that ONE must
   now choose meager sets from ${\cal A}$ only. The proof of
   Theorem \ref{countablefin-neg} gives a 
   cofinal family ${\cal A}$ such that TWO does not have a winning
   remainder strategy in the game $WMG({\cal A},[\omega_1]^{<\aleph_0})$.
   This should be contrasted with Theorem \ref{wmeg-pos2}, which implies
   that there are many uncountable cardinals $\kappa$ such that for some
   cofinal family ${\cal A}\subset[\kappa]^{\aleph_0}$, TWO has a winning
   remainder strategy in $WMG({\cal A},[\kappa]^{<\aleph_0})$.

\begin{theorem}\label{codingtoremainder}If TWO has a winning coding
   strategy in $WMG(J)$, and if there is a 
   cofinal family ${\cal A}\subset\langle J\rangle$ such that $A\Delta
   B\not\in J$ whenever $A\neq B$ are elements of ${\cal A}$, then TWO
   has a winning remainder strategy in $WMG({\cal A},J)$.
\end{theorem}

\begin{description}\item[Proof]{Let $F$ be a winning coding strategy
   for TWO in $WMG(J)$, and let ${\cal A}\subset\langle J\rangle$ be a
   cofinal family as in the hypothesis of the theorem. If $B$ is not in
   ${\cal A}$, but there is an $A\in {\cal A}$ such that $B\subset A$ and
   $A\backslash B\in J$, then this $A$ is unique on account of the
   properties of ${\cal A}$.

   Define a remainder strategy $G$ for TWO as follows: Let $B\in\langle
   J\rangle$ be given.
\begin{enumerate}
\item{$G(B)=F(\emptyset,B)$ if $B\in{\cal A}$,}
\item{$G(B)=F(A\backslash B,A)$ if $B\not\in{\cal A}$, but $B\subset
   A$ and $A\backslash B\in J$ for an $A\in {\cal A}$, and}
\item{$G(B)=\emptyset$ in all other circumstances.}
\end{enumerate}
   Then $G$ is a winning remainder strategy for TWO in $WMG({\cal A},J)$.
   $\QED$}
\end{description}

\begin{corollary}\label{smg1}Let $\lambda$ be a cardinal number of
   countable cofinality. For each $\kappa\geq\lambda$, there is a cofinal
   family ${\cal A}\subset[\kappa]^{\lambda}$ such that TWO has a winning
   remainder strategy in $WMG({\cal A},J)$.
\end{corollary}

\begin{description}\item[Proof]{ Write $\kappa=\cup_{\alpha<\kappa}
   S_{\alpha}$ where $\{S_{\alpha}:\alpha<\kappa\}$ is a disjoint
   collection of sets, each of cardinality $\lambda$. For each
   $A\in[\kappa]^{\lambda}$, put $A^*=\cup_{\alpha\in A}S_{\alpha}$. Then
   ${\cal A}=\{A^*:A\in[\kappa]^{\lambda}\}$ is a cofinal subset of
   $[\kappa]^{<\lambda}$ which has the properties required in Theorem 
   \ref{codingtoremainder}. The result now follows from that theorem and
   the fact that TWO has a winning coding strategy in
   $WMG([\kappa]^{<\lambda})$ - see \cite{S4}. $\QED$}
\end{description}

   It is worth noting that Corollary \ref{smg1} is a result in
   ordinary set theory, whereas we used the Generalized Continuum
   Hypothesis in Corollary \ref{wmegcor}. 

\section{The very strong game, VSG(J).}

   Moves by player TWO in the game $VSG(J)$ (introduced in
   \cite{B-J-S}) consist of pairs of the
   form $(S,T)\in\langle J\rangle\times J$, while those of ONE are
   elements of $\langle J\rangle$. A sequence
$$(O_1,(S_1,T_1), O_2,(S_2,T_2),\dots)$$ 
   is a play of $VSG(J)$ if:
\begin{enumerate}
\item{$O_{n+1}\supseteq S_n\cup T_n$, and }
\item{$O_n, S_n\in\langle J\rangle$ and $T_n\in J$ for each $n$.}
\end{enumerate}
   Player TWO wins such a play if
$$\cup_{n=1}^{\infty}O_n\subseteq\cup_{n=1}^{\infty}T_n.$$
   A strategy $F$ is a remainder strategy for TWO in $VSG(J)$ if
$$(S_{n+1},T_{n+1})=F(O_{n+1}\backslash(\cup_{j=1}^{n}T_n))$$
   for each $n$.   

   For $X\in\langle J\rangle$ we write $F(X)=(F_1(X),F_2(X))$ when $F$
   is a remainder strategy for TWO in $VSG(J)$. When $F$ is a winning
remainder strategy for TWO, we may assume that it has the following
properties:
\begin{enumerate}
\item{$F_1(X)\cap F_2(X)=\emptyset$; for $G$ is a
   winning remainder strategy if $G_1(X)=F_1(X)\backslash F_2(X)$ and
   $G_2(X)=F_2(X)$ for each $X$.}
\item{$X\backslash F_2(X)\subseteq F_1(X)$; for $G$ is a winning
   remainder strategy if $G_1(X)=(X\cup F_1(X))\backslash F_2(X)$ and
   $G_2(X)=F_2(X)$ for each $X$.}
\end{enumerate}
   The following Lemma describes a property which every winning
   remainder strategy of player TWO for the game $VSG(J)$ must have.

\begin{lemma}\label{property1}Assume that $J\subset\langle
   J\rangle\subset{\EuScript P}(S)$ and let $F$ be a
   winning remainder strategy for TWO in the game $VSG(J)$. Then the
   following assertion holds.
\begin{quote} For each $x\in S$ there exist a $C_x\in\langle J\rangle$
   and a $D_x\in J$ such that:
\begin{enumerate}
\item{$C_x\cap D_x=\emptyset$ and}
\item{$x\in F_2(B)$ for each $B\in\langle J\rangle$ such that
   $C_x\subseteq B$ and $D_x\cap B=\emptyset$.}
\end{enumerate}
\end{quote}
\end{lemma}

\begin{description}\item[Proof]{Let $F$ be a remainder strategy of
   TWO, but assume the negation of the conclusion of the lemma. We
   also assume that for each $X\in\langle J\rangle$, $X\backslash
   F_2(X)\subseteq F_1(X)$ and $F_1(X)\cap F_2(X)=\emptyset.$

   Choose an $x\in S$ witnessing this negation. Then there is for each
   $C\in\langle J\rangle$ and for each $D\in J$ with $x\in C$ and $C\cap
   D=\emptyset$ a $B\in\langle J\rangle$ such that $B\cap D=\emptyset$,
   $C\subseteq B$ and $x\not\in F_2(B)$. We now construct a
   sequence $\langle(B_k,C_k,D_k,M_k,S_k,N_k):k\in\naturals\rangle$ as
   follows: (we go through the first three steps of the construction
   for clarity, before stating the general requirements for the sequence)

   Put $C_1=\{x\}$ and $D_1=\emptyset$. Choose $B_1\in\langle
   J\rangle$ such that $C_1\subseteq B_1$ and $x\not\in F_2(B_1)$.
   Put $M_1=B_1$ and $(S_1,N_1)=F(M_1)$.
   This defines $(B_1,C_1,D_1,M_1,S_1,N_1)$.

   Put $C_2=S_1$ and $D_2=N_1$. Choose $B_2\in\langle J\rangle$ such
   that $C_2\subseteq B_2$, $D_2\cap B_2=\emptyset$, and $x\not\in
   F_2(B_2)$. Put $M_2=B_2\cup D_2$ and $(S_2,N_2)=F(M_2\backslash N_1)$.
   This defines $(B_2,C_2,D_2,M_2,S_2,N_2)$.

   Put  $D_3=(N_1\cup N_2)$ and $C_3=S_2\backslash D_3$. Choose
   $B_3\in\langle J\rangle$ such that $C_3\subseteq B_3$, $D_3\cap
   B_3=\emptyset$, and $x\not\in F_2(B_3)$. Put $M_3=B_3\cup D_3$ and 
   $(S_3,N_3)=F(M_3\backslash D_3)$. This defines
   $(B_3,C_3,D_3,M_3,S_3,N_3)$. 

   Let $k\geq 3$ be given, and assume that 
\[(B_1,C_1,D_1,M_1,S_1,N_1),\dots,(B_k,C_k,D_k,M_k,S_k,N_k)\] 
   have been chosen so that:
\begin{enumerate}
\item{$D_{j+1}=(N_1\cup\dots\cup N_j)\in J$,}
\item{$C_{j+1}=S_{j}\backslash D_{j+1}\in\langle J\rangle$ and
   $x\in C_{j+1}$ for $j<k$,}
\item{$C_j\subseteq B_j$, while $x\not\in F_2(B_j)$ and $B_j\cap
   D_j=\emptyset$, and}
\item{$M_{j}=B_{j}\cup D_{j}$ and }
\item{$(S_j,N_j)=F(M_j\backslash D_j)$ for $j\leq k$, and}  
\item{$(B_1,C_1,D_1,M_1,S_1,N_1)$ and $(B_2,C_2,D_2,M_2,S_2,N_2)$ are
   as above.} 
\end{enumerate}

   Put $D_{k+1}=\cup_{j=1}^kN_j$ and $C_{k+1}=S_k\backslash D_{k+1}$. By
   hypothesis there is a $B_{k+1}\in\langle J\rangle$ such that
   $C_{k+1}\subseteq B_{k+1}$, $B_{k+1}\cap D_{k+1}=\emptyset$, and
   $x\not\in F_2(B_{k+1})$. Put $M_{k+1}=B_{k+1}\cup D_{k+1}$, and put
   $(S_{k+1},N_{k+1})=F(M_{k+1}\backslash(N_1\cup\dots\cup N_k))$.

   Then the sequence
\[(B_1,C_1,D_1,M_1,S_1,N_1), \dots,
   (B_{k+1},C_{k+1},D_{k+1},M_{k+1},S_{k+1},N_{k+1})\] 
   still satisfies properties 1 through 6 above. Continuing like this we
   obtain an infinite sequence 
   which satisfies these. But then 
\[(M_1,(S_1,N_1),\dots,M_k,(S_k,N_k),\dots)\]
   is a play of $VSG(J)$ during which player TWO used the
   remainder strategy $F$; this play is moreover lost by TWO because TWO
   never covered the point $x$. This proves the contrapositive of the
   lemma. $\QED$}
\end{description}

\begin{theorem}[Galvin]\label{countablefin2} For $\kappa>\aleph_1$,
   TWO does not have a winning remainder strategy in
   $VSG([\kappa]^{<\aleph_0})$.
\end{theorem}

\begin{description}\item[Proof]{ Let $F$ be a remainder strategy for
   TWO. If it were winning, choose for each $x\in\kappa$ a
   $D_x\in[\kappa]^{<\aleph_0}$ and a $C_x\in[\kappa]^{\leq \aleph_0}$
   such that:
\begin{enumerate}
\item{$C_x\cap D_x=\emptyset$,}
\item{$x\in C_x$ and}
\item{$x\in F_2(B)$ for each $B\in[\kappa]^{\leq\aleph_0}$ such that $B\cap
   D_x=\emptyset$ and $C_x\subseteq B$.}
\end{enumerate}

   Now $(D_x:x\in\kappa)$ is a family of finite sets. By the
   $\Delta$-system lemma we find an $S\in[\kappa]^{\kappa}$ and a finite
   set $R$ such that $(D_x:x\in S)$ is a $\Delta$-system with root $R$.
   For $x\in S$ define:
\[f(x)=\{y\in S:D_y\cap C_x\neq\emptyset\}.\]
   Then $f(x)$ is a countable set and $x\not\in f(x)$ for
   each $x\in S$. By Hajnal's set-mapping theorem we find
   $T\in[S]^{\kappa}$ such that $C_x\cap D_y=\emptyset$ for all $x,y\in
   T$. 

   Let $K\in[T]^{\aleph_0}$ be given, and put $B=\cup_{x\in K}C_x$. Then
   $K\subseteq F_2(B)$, a contradiction. $\QED$}
\end{description}

   Using similar ideas but with the appropriate cardinality assumption to
   ensure that the corresponding versions of the $\Delta$-system lemma
   and the set-mapping theorems are true, one obtains also:

\begin{theorem}\label{uncountablecofneg} Let $\lambda$ be a cardinal
   of countable cofinality. If
   $\kappa>2^{\lambda}$, then TWO does not have a winning remainder
   strategy in $\vsgkl$. 
\end{theorem}

    Since for every cardinal $\lambda$ of
   countable cofinality, and for each cardinal $\kappa$ player TWO has a
   winning coding strategy in $WMG([\kappa]^{<\lambda})$ (see for
   example \cite{S4}), Theorems \ref{countablefin2} and \ref{uncountablecofneg}
   also show that the existence of a winning remainder strategy
   for TWO in $VSG(J)$ is a stronger statement than the existence of a
   winning coding strategy for TWO.

\begin{problem}Let $\lambda$ be an uncountable cardinal of countable
   cofinality. Let $\kappa$ be an infinite cardinal number such that
   $\lambda^{<\lambda}<cof([\kappa]^{\lambda},\subset) \leq
   2^{\lambda}.$
   Does TWO fail to have a winning remainder strategy in any of 
   $WMEG([\kappa]^{<\lambda})$, $WMG([\kappa]^{<\lambda})$ or
   $VSG([\kappa]^{<\lambda})$? 
\end{problem}

   The following theorem shows that the $\kappa$ in Theorem
   \ref{countablefin2} cannot be decreased to $\omega_1$. Thus,
   the rules of the very strong game are more advantageous to TWO
   than those of the other versions we considered earlier in this paper.

\begin{theorem}\label{vsg2} If $cof(\langle
   J\rangle,\subset)=\aleph_1$, then TWO has a winning 
   remainder strategy in $VSG(J)$.
\end{theorem}

\begin{description}\item[Proof]{ We may assume that there is for each
   $X\in\langle J\rangle\backslash J$, a $Y\in\langle J\rangle\backslash
   J$ such that $X\cap Y=\emptyset$ (else, TWO has an easy winning
   remainder strategy even in $WMEG(J)$). Let $\prec$ be a well-ordering
   of $S$, the underlying set of our topological space. Choose two
   $\omega_1$-sequences $(C_{\alpha}:\alpha<\omega_1)$ and
   $(x_{\alpha}:\alpha<\omega_1)$ such that:
\begin{enumerate}
\item{$C_{\alpha}\subset C_{\beta}\in\langle J\rangle$,}
\item{$x_{\alpha}\in C_{\beta}$,}
\item{$x_{\beta}\not\in C_{\beta}$,}
\item{$x_{\alpha}\prec x_{\beta}$ and}
\item{$C_{\beta}\backslash C_{\alpha}\not\in J$}
   for all $\alpha<\beta<\omega_1$, and
\item{$\{C_{\alpha}:\alpha<\omega_1\}$ is cofinal in $\langle
   J\rangle$.}
\end{enumerate}

   For each $X\in\langle J\rangle$ we write $\beta(X)$ for
   $\min\{\alpha<\omega_1:X\subseteq C_{\alpha}\}$. Put ${\frak
   X}=\{x_{\alpha}:\alpha<\omega_1\}$. Write $\Omega$ for
   $\omega_1\backslash\omega$. Let $F$ be a winning perfect
   information strategy for TWO in $REG(J)$, and let $G$ be a
   winning perfect information strategy for TWO in
   $REG([\{x_{\delta}:\delta\in\Omega\}]^{<\aleph_0})$. We may assume
   that if $\sigma$ is a 
   sequence of length $r$ of subsets of $\Omega$, at least one of which
   is infinite, then $|G(\sigma)|\geq r$. We also define:
   $K_{\beta}=\{x_{\gamma}:\gamma\leq\beta\}$ for each $\beta\in\Omega$.

   We define a remainder strategy $H$ for TWO in $VSG(J)$. Thus, let
   $B\in\langle J\rangle$ be given.
\begin{enumerate}
\item{If $B\in J$: Then put $H(B)=(C_{\beta(B)+\omega},\{x_0,x_{\beta(B)}\})$}
\item{If $B\not\in J$:
\begin{enumerate}
\item{If $\{n<\omega:x_n\not\in B\}=\{0,1,\dots,k\}$:\\
   Let $T$ be $\{x_{\beta(B)}\}$ together with the first $\leq k+1$
   elements of $\{x_{\alpha}:\alpha\in\Omega\}\backslash B$. Put
\[S=T\cup(\cup\{G(\sigma):\sigma\in\mbox{$^{\leq
   k+2}(\{K_{\delta}:x_{\delta}\in T$}\}),
\]
   a set in $[\{x_{\delta}:\delta\in\Omega\}]^{<\aleph_0}$. Let
   $p$ be the cardinality of $S$. Then define
\[\overline{S}=\{x_0,\dots,x_p\}\cup S \cup((\cup\{F(\sigma):\sigma\in
   \mbox{$^{\leq p}(\{C_{\alpha}:x_{\alpha}\in S\})$}\})\backslash{\frak X}).
\]
   Put $H(B)=(C_{\beta(B)+\omega},\overline{S})$.}
\item{If $\{n<\omega:x_n\not\in B\}$ is not a finite initial segment
   of $\omega$: Then we put $H(B)=(C_{\beta(B)+\omega},\{x_0,x_{\beta(B)}$.}
\end{enumerate}}
\end{enumerate}

   To see that $H$ is a winning remainder strategy for TWO, consider a
   play
$$(O_1,(S_1,T_1),\dots,O_n,(S_n,T_n),\dots)$$ 
   where $(S_1,T_1)=H(O_1)$ and
   $(S_{n+1},T_{n+1})=H(O_{n+1}\backslash(\cup_{j=1}^n S_j)$ for each $n$.

   For convenience we put 
\begin{itemize}
\item{$W_0=T_0=\emptyset$ and $W_{n+1}=W_n\cup T_{n+1}$,} 
\item{$B_n=O_n\backslash W_n$, }
\item{$\beta_n=\beta(B_n)$ and $\alpha_n=\beta_n+\omega$}
\end{itemize}
   for each $n$

   Note that if $B_j$ is such that $\{n\in\omega:x_n\not\in
   B_j\}(=\{0,1,\dots,k_j\}$ say) is a finite initial segment of
   $\omega$, then the same is true for $B_{j+1}$. It follows that
   $(S_j,T_j)$ is defined by Case 2(a) for each $j>1$, and that
   $(k_j:j\in\naturals)$ is an increasing sequence. It further follows
   that $\{x_{\beta_1},\dots,x_{\beta_j}\}\subset T_j$ for these $j$.
   This in turn implies that:
\begin{enumerate}
\item{$\cup_{j=1}^{\infty} F(C_{\beta_1},\dots,C_{\beta_j})\backslash{\frak
   X}\subseteq \cup_{n=1}^{\infty}T_n$,}
\item{$\cup_{j=1}^{\infty} G(K_{\beta_1},\dots,K_{\beta_j})\subseteq
   \cup_{n=1}^{\infty}T_n$, and}
\item{$\{x_n:n<\omega\}\subseteq\cup_{n=1}^{\infty}T_n$.}
\end{enumerate}

   But then $\cup_{n=1}^{\infty}O_n\subseteq\cup_{n=1}^{\infty}T_n$.
$\QED$}
\end{description}

\begin{corollary}\label{vsgom1}TWO has a winning remainder strategy
   in $VSG([\omega_1]^{<\aleph_0})$
\end{corollary}

   Using the methods of this paper we can also show that if
   $J\subset{\EuScript P}(S)$ is a free ideal such that
   there is an $A\in\langle J\rangle$ such that
   $cof(\langle J\rangle,\subset)\leq |J\lceil_A|$,
   then TWO has a winning remainder strategy in $VSG(J)$.

\begin{corollary}\label{vsgom2} For every $T_1$-topology on
   $\omega_1$, without isolated points, TWO has a winning remainder
   strategy in $VSG(J)$.
\end{corollary}

\end{document}